# About a virtual subset

Dipl. Math.(FH) Klaus Lange*

*Cinterion Wireless Modules GmbH, Technology, Product Development, Integration Test
Berlin, Germany
Klaus.Lange@cinterion.com

___

Abstract: Two constructed prime number subsets (called "prime brother & sisters" and "prime cousins") lead to a third one (called "isolated primes") so that all three disjoint subsets together generate the prime number set. It should be suggested how the subset of isolated primes give a new approach to expand the set theory by using virtual subsets.

___

## I.    Prime number brothers & sisters

This set of prime numbers is given

**P** = {2; 3; 5; 7; 11; …}.

Firstly we are going to establish a subset of prime numbers called "brothers & sisters" to generalize the well known prime number twins.

Definition D1: Brother & Sister Primes

Two direct neighbours of prime numbers $p_i$ and $p_{i+1}$ are called brothers & sisters if this distance d is given

$d = p_{i+1} - p_i = 2^n$ with $n \in \mathbb{IN}_0$.

All these brothers & sisters are elements of the set

**B** := $\cup B_j$

Notes:
(i) $B_j$ is the j-th subset of **B** by d and separated from $B_{j+1}$ because the distance between $B_j$ and $B_{j+1}$ does not have the structure $2^n$. For example

**B**$_1$ = {2; 3; 5; 7; 11; 13; 17; 19; 23}
**B**$_2$ = {29; 31}
**B**$_3$ = {37; 41; 43; 47}
**B**$_4$ = {59; 61}
**B**$_5$ = {67; 71; 73}
**B**$_6$ = {79; 83}
**B**$_7$ = {89; 97; 101; 103; 107; 109; 113}
**B**$_8$ = {127; 131}
**B**$_9$ = {137; 139}
**B**$_{10}$ = {149; 151}

The distances between 29 - 23 and 37 – 31 etc. do not fulfil the form $d = 2^n$. That is the reason why prime number 53 is not listed.

(ii) The prime brothers & sisters are a general structure of the twin prime distance, therefore we can assume that infinite many subsets $B_j$ exist. But it has not been proven yet, we will discuss this topic later.

(iii) It is an interesting question how many elements can have a subset $B_j$. Later on we will proof that there is no last $B_j$ which has infinite many elements.

## II.  Prime number cousins

Now we are going to establish a disjoint prime number set to **B**, so that the combination of this new subset of "other" primes **O** and the subset **B** rebuild the set **P**.

Definition D2: Other Primes

**O** : = **P/B**.
**O** are called the "Other" primes in **P** without **B**.

Notes:
(iv) Following D2 it is **O** ∩ **B** = { }.

(v) We remark that **O** has the same structure of subsets $O_k$ like **B** which has the subsets $B_j$. For example the first ten subsets of **O** are these:

**O**$_1$ = {53}
**O**$_2$ = {157}
**O**$_3$ = {173}
**O**$_4$ = {211}
**O**$_5$ = {251; 257; 263}

$O_6 = \{293\}$
$O_7 = \{331; 337\}$
$O_8 = \{373\}$
$O_9 = \{509\}$
$O_{10} = \{541; 547; 557; 563\}$

(vi) To give a first idea how the subsets of **O** fill the gaps of **B**, so that all prime numbers can be listed, we have a short look to all prime numbers up to number 300 using both subsets of **P**:

$B_1 = \{2; 3; 5; 7; 11; 13; 17; 19; 23\}$
$B_2 = \{29; 31\}$
$B_3 = \{37; 41; 43; 47\}$
$O_1 = \{53\}$
$B_4 = \{59; 61\}$
$B_5 = \{67; 71; 73\}$
$B_6 = \{79; 83\}$
$B_7 = \{89; 97; 101; 103; 107; 109; 113\}$
$B_8 = \{127; 131\}$
$B_9 = \{137; 139\}$
$B_{10} = \{149; 151\}$
$O_2 = \{157\}$
$B_{11} = \{163; 167\}$
$O_3 = \{173\}$
$B_{12} = \{179; 181\}$
$B_{13} = \{191; 193; 197; 199\}$
$O_4 = \{211\}$
$B_{14} = \{223; 227; 229; 233\}$
$B_{15} = \{239; 241\}$
$O_5 = \{251; 257; 263\}$
$B_{16} = \{269; 271\}$
$B_{17} = \{277; 281; 283\}$
$O_6 = \{293\}$

That indirect definition is needed to establish positive properties for a counter subset in contrast to **B** by using the same distance structure $2^n$.

Definition D3: Cousin Primes

Given primes p, q ∈ **O** with p > q . This primes are called Cousins if they fulfil the distance structure

$d = p - q = 2^n$.

**C** := $\{p > q \mid p; q \in \mathbf{O} \wedge p - q = 2^n\}$

Notes:
(vi) The existence of **C** is shown by some examples of $O_k$

$O_2 = \{157\}$

$O_3 = \{173\}$
$O_9 = \{509\}$
$O_{10} = \{541; 547; 557; 563\}$

Because of

$173 - 157 = 16 = 2^4$
$541 - 509 = 32 = 2^5$
$557 - 541 = 16 = 2^4$
$563 - 547 = 16 = 2^4$

we can write temporary

$C_i = \{157; 173; …\}$
$C_{i+1} = \{509; 541; 557; …\}$
$C_{i+2} = \{547; 563; …\}$

(vii) Because of **O** it is valid that $C \cap B = \{\}$. In other words: Prime numbers which have brothers and sisters do not have cousins and vice versa.

(viii) The next question: **Is $P = C \cup B$ true?**

## III. Discussion about the existence of isolated prime numbers

Discussing the question of (viii) we can set up the distance relationship between prime numbers in general.

Definition D4: Relative primes

Every prime p which has a distance $d = 2^n$ to another prime q is member of Relative prime set **R**. The numbers p and q called Relative primes.

Notes:
(ix) Following D4 not only neighbour primes with distance $2^n$ are members of **R**. All primes which have one partner in d at least are elements of **R**. It is clear that Relative primes reunion the subsets **B** and **C**. **R** implies the same question as (viii) just in another form:
**Is $P = R$ ?**

(x) The same question just in another form highlights the main idea of the construction of subset **C** in contrast to set **R**:
**Is $R = C \cup B$ ?**

In deep discussion:

First of all the existence of a prime number p which has no partner prime q for

$|p - q| = 2^n$

is theoretical possible.

The same number p which is isolated by D4 will also be isolated by D3 in the same way. Following D4 it will be difficult to detect some candidates for isolated prime numbers, using D3 it is much easier.

For example we have a look to the prime number 53 given in $\mathbf{O}_1$:

D4 allows a distance relation independent of $\mathbf{B}$, so the prime numbers in $\mathbf{B}_3$ and $\mathbf{B}_4$ are useable without limitations.

$53 - 37 = 16 = 2^4$

or

$61 - 53 = 8 = 2^3$

etc.

$\qquad$ D4 $\Rightarrow$ 53 $\in$ $\mathbf{R}$.

Following D3 it is forbidden to build the distance from 53 to a partner prime in $\mathbf{B}$.

Searching for a partner prime under limitations of D3 - only in $\mathbf{C}$ - is very difficult.

All p < 53 are elements of $\mathbf{B}$, so they are not useable.

Reaching out for p > 53 this results are shown:

53 +…
$2^1$ = 55 no prime
$2^2$ = 57 no prime
$2^3$ = 61 $\in$ $\mathbf{B}$ (brother: 59)
$2^4$ = 69 no prime
$2^5$ = 85 no prime
$2^6$ = 117 no prime
$2^7$ = 181 $\in$ $\mathbf{B}$ (brother: 179)
$2^8$ = 309 no prime
$2^9$ = 565 no prime
$2^{10}$ = 1077 no prime
$2^{11}$ = 2101 no prime
$2^{12}$ = 4149 no prime
$2^{13}$ = 8245 no prime
$2^{14}$ = 16437 no prime
$2^{15}$ = 32821 no prime
$2^{16}$ = 65589 no prime
$2^{17}$ = 131125 no prime
$2^{18}$ = 262197 no prime
$2^{19}$ = 524341 $\in$ $\mathbf{B}$ (brother: 524309)
$2^{20}$ = 1048629 no prime
$2^{21}$ = 2097205 no prime
$2^{22}$ = 4194357 no prime
$2^{23}$ = 8388661 no prime

$2^{24}$ = 16777269 no prime
$2^{25}$ = 33554485 no prime
$2^{26}$ = 67108917 no prime
$2^{27}$ = 134217781 ∈ **B** (brother: 134217779)
$2^{28}$ = 268435509 no prime
$2^{29}$ = 536870965 no prime
$2^{30}$ = 1073741877 no prime
$2^{31}$ = 2147483701 no prime
$2^{32}$ = 4294967349 no prime
$2^{33}$ = 8589934645 no prime
$2^{34}$ = 17179869237 no prime
$2^{35}$ = 34359738421 prime number without a Brother

Result: 53 and 53 + $2^{35}$ are elements of the Cousin prime number subset.

If there is no partner prime in **O** which has the distance $2^n$ we call it an Isolated prime number.

Definition D5 Isolated Prime subset

Every prime number which is no element of **B** or **C** is called an Isolated prime number and is a member of the set **I**.

**I** := {p | p ∈ **P**/(**B**∪**C**) }.

Notes:
(xi) It is possible that p ∈ **I** and q ∈ **B** with d = |q - p| = $2^n$.

(xii) Following D5, D3 and D1 it is obvious that **P** = **B** ∩ **C** ∩ **I**.

(xiii) Proving p ∈ **I** we need many infinite steps of calculation searching for |q - p| ≠ $2^n$ with q ∈ **O** = **P**/**B**.

Note (xiii) leads up to a need for a criterion of only finite steps of calculation. Finding this, we are able to call a prime number as a candidate of **I**.

Definition D6 Isolated prime candidate

The prime number p ∈ **P**/**B** is called a candidate of **I** if q ∈ **P**/**B**; q < p with p – q ≠ $2^n$ and all sums p + $2^i$ ∉ **P** with i = 1 to p.

Example:

Testing the prime number 211 there is no prime number q ∈ **P**/**B** with p – q = $2^n$. And all sums from 211 + 2 up to 211 + $2^{211}$ are either no primes or primes in **B**. The prime number 211 is a candidate for the subset **I**. Going on with brute force computing we will find the sum 211 + $2^{448}$ ∈ **P**/**B**. That proves p = 211 and q = 211 + $2^{448}$ ∈ **C**.

# IV. Analysing by combination the subsets B, C and I

The prime subsets **B**, **C** and **I** are useable to analyse **P**.

So we give the

Definition D7 Set criteria function

Be **X** a given set then we define a function

$$\psi(\mathbf{X}) = \begin{cases} -1, & \mathbf{X} = \{\} \\ 0, & \text{for } \mathbf{X} \text{ has a finite number of elements} \\ 1, & \text{for } \mathbf{X} \text{ has an infinite number of elements} \end{cases}$$

and

Definition D8 Possible combination function

Be the values of $\psi(\mathbf{X}_i)$, $i = 1, 2, 3$ are given then we define a function

$$\kappa_j(\psi(\mathbf{X}_1), \psi(\mathbf{X}_2), \psi(\mathbf{X}_3)) = \begin{cases} 0, & \text{if the combination of } \psi \text{ - values is not possible} \\ 1, & \text{other} \end{cases}$$

To do:
Our task is to find those combinations of $\psi(\mathbf{X}_i)$ with

$\kappa(\psi(\mathbf{X}_1), \psi(\mathbf{X}_2), \psi(\mathbf{X}_3)) = 1$.

The examples of elements in **B** and **C** show, that

$-1 < \psi(\mathbf{B}),\ \psi(\mathbf{C}) \leq 1$.

Only **I** fits to all three possible values of function $\psi$

$-1 \leq \psi(\mathbf{I}) \leq 1$.

So we have twelve combinations to analyse:

$\kappa_j = \kappa(\psi(\mathbf{B}), \psi(\mathbf{C}), \psi(\mathbf{I}))$ with

$$\kappa_1 = \kappa(0, 0, -1); \kappa_2 = \kappa(0, 0, 0); \kappa_3 = \kappa(0, 0, 1)$$
$$\kappa_4 = \kappa(0, 1, -1); \kappa_5 = \kappa(0, 1, 0); \kappa_6 = \kappa(0, 1, 1)$$
$$\kappa_7 = \kappa(1, 0, -1); \kappa_8 = \kappa(1, 0, 0); \kappa_9 = \kappa(1, 0, 1)$$
$$\kappa_{10} = \kappa(1, 1, -1); \kappa_{11} = \kappa(1, 1, 0); \kappa_{12} = \kappa(1, 1, 1).$$

Statement S1:

$\kappa_1 = \kappa_2 = \kappa_3 = \kappa_5 = \kappa_6 = 0.$

Proof Arguments PA1:

In the case of $\kappa_1 = \kappa_2 = 0$ is the argument clear, because for prime number set **P** we have

$\psi(\mathbf{P}) = 1$ and $\mathbf{P} = \mathbf{B} \cap \mathbf{C} \cap \mathbf{I}$

it follows

$\psi(\mathbf{B} \cap \mathbf{C} \cap \mathbf{I}) = \psi(\mathbf{B}) \cup \psi(\mathbf{C}) \cup \psi(\mathbf{I}) = 1$

in contradiction to

$\kappa_1 : \psi(\mathbf{B}) \cup \psi(\mathbf{C}) \cup \psi(\mathbf{I}) = 0 \cup 0 \cup -1 < 1.$
$\kappa_2 : \psi(\mathbf{B}) \cup \psi(\mathbf{C}) \cup \psi(\mathbf{I}) = 0 \cup 0 \cup 0 < 1.$

In the case of $\kappa_3 = \kappa_5 = \kappa_6 = 0$ we argue heuristically. Because of the prime number theorem the possibility $\varphi$ of the number $x \in \{1, 2, 3, \ldots, n\}$ is a prime number given by

$\varphi_{(x \text{ is prime})} \approx 1/(\ln(n))$

$\ln(n)$ is the natural logarithm of n.

Be $p \in \mathbf{P}$. If m large enough we have

$q = p + 2^m < 2^{m+1}$

it follows

$\varphi_{(q \in \mathbf{P})} > 1/((m+1)\ln(2)).$

The summation over all these large enough $m \to \infty$ leads to

$$S = \sum_{m \to \infty} 1/((m+1)\ln(2))$$

and that is a harmonic series.

Because of the divergence of S it is clear:
It exists not only one but infinite many $q = p+2^m \in \mathbf{P}$.

That means: If we do not want to contradict the Definition D6 we need infinite many $q \in \mathbf{B}$ for only one $p \in \mathbf{I}$. For $\kappa_j$ with values for $\psi(\mathbf{I}) > -1$ it must be $\psi(\mathbf{B}) = 1$ following D6. That is the reason why $\kappa_3 = \kappa_5 = \kappa_6 = 0$ is heuristically correct.

Statement S2:

$\kappa_7 = \kappa_8 = \kappa_9 = 0$.

Proof Arguments PA2:

Case $\kappa_7 = \kappa_8 = 0$:
$\psi(\mathbf{B}) = 1$, $\psi(\mathbf{C}) = 0$ and $\psi(\mathbf{I}) < 1$ means, that there is a $\varepsilon \gg 0$ and for $p > \varepsilon$ it is $p \in \mathbf{B}$.
Using the theorem of P.G.L. Dirichlet we have in statistical approximation the same numbers of primes in both forms $4n + 1$ and $4n + 3$ if n is large enough. So we can use $n \gg \varepsilon$ to fulfil this condition and for those n all primes $p = 4n \pm 1 \in \mathbf{B}$.
There are only two ways for $p > \varepsilon$ to change from $4n+1$ to $4n-1$:

a) Many tiny subsets $\mathbf{B}_j$ needing the distance $\delta_B = 4n+2$ from the last prime of $\mathbf{B}_j$ to the first prime of $\mathbf{B}_{j+1}$. Statistically we have the same number of $\delta_B = 4n+2$ and $\delta_B = 4n$. That means: Changing forms far enough using distances between subsets of $\mathbf{B}$ is not enough to fulfil Dirichlet theorem.
b) To satisfy Dirichlet theorem we need changing forms in the midst of $\mathbf{B}_j$. For this only the twin prime gaps $d = p_{i+1} - p_i = 2$ help.

Using the twin prime constant $\beta_2$ for the heuristically analysis

$$\beta_2 = 2\prod_{p>2}[1 - (p-1)^{-2}] \approx 1{,}320$$

we see that the twin primes for very large $p \gg \varepsilon$ are very rare and the statistically approximated equilibrium of $4n + 1$ and $4n - 1$ is strongly disrupted (not only weak in the sense of Chebyshev's bias) if there are only infinitely many primes in $\mathbf{B}$, but not in $\mathbf{C}$ and not in $\mathbf{I}$.

Case $\kappa_9 = 0$:
In that case we have infinitely many primes in $\mathbf{B}$ and in $\mathbf{I}$. But in PA1 we saw that

$\pi_\mathbf{B}(n) \gg \pi_\mathbf{I}(n)$

for $\pi_\mathbf{X}(n)$ is the number of elements up to n for set X:
Every $i \in \mathbf{I}$ needs infinity many $b \in \mathbf{B}$ and for $\psi(\mathbf{I}) = 1$ the disruption of statistically approximated equilibrium regarding $4n + 1$ and $4n - 1$ is much stronger than in the case before, for reasons of the twin prime constant never changes and the rare i's cannot help to heal this stronger disruption.

As result of S1 and S2, we have

Statement S3:

$\psi(\mathbf{C}) = 1$.

As a first result we see that only four combinations are possible:

$\kappa_4$; $\kappa_{10}$; $\kappa_{11}$; $\kappa_{12}$

$\kappa_4$ is the only one of this four combinations which has $\psi(\mathbf{B}) = 0$ and $\kappa_{12}$ is the only one which has $\psi(\mathbf{I}) = 1$.

That leads to the conjecture $\psi(\mathbf{B}) = 1$ and $\psi(\mathbf{I}) < 1$.

Under precondition of the unproven Elliott – Halberstam conjecture Goldston, Pintz and Yildirim show that

$\lim_{n \to \infty} d_n \leq 16 = 2^4$

which is known as the bounded gap conjecture and for $d = 2^4$ it confirms

Statement S4:

Under bounded gap conjecture it is

$\psi(\mathbf{B}) = 1$.

As we saw before the existence of the Isolated prime number set is improvable in the sense of Gödel theorem. But we can see that $\psi(\mathbf{I}) < 1$ is given.

Be $\psi(\mathbf{I}) = 1$ then we have infinite many elements from $\mathbf{I}$. All distances $d = 2^n$ for every $i \in \mathbf{I}$ with $j = i \pm 2^n$ point at a number which is either not a prime number or a prime number in $\mathbf{B}$. It is not allowed that $j \in \mathbf{I}$ or $j \in \mathbf{C}$. To reduce the probability to zero that j is in $\mathbf{C}$ or $\mathbf{I}$, there must be infinite elements in $\mathbf{I}$ and $\mathbf{C}$ which are very very rare for large n's in $j = i + 2^n$. In other words:

If n is large enough then it must be fulfilled

$\pi_\mathbf{B}(n) \gg \pi_\mathbf{I}(n)$

and

$\pi_\mathbf{B}(n) \gg \pi_\mathbf{C}(n)$.

Looking at this we have the same problem to satisfy Dirichlet theorem about the numbers of prime forms 4n+1 and 4n+3 as in PA2. We conclude:

Statement S5:

$\psi(\mathbf{I}) < 1$.

Following Gödel we are not able to decide which of the combination $\kappa_{10}$ and $\kappa_{11}$ is true.

## V. A virtual subset

The constructed subset **I** is an example of a virtual subset.

My suggestion to define those subsets is:

Definition D9 Virtual subset

A virtual subset $V \subset W$, $W \neq \{\}$ is given by the following criterions

a) $W \neq V$
b) $W = U_1 \cup U_2 \cup \ldots \cup U_n \cup V$ with $U_i \cap U_j = \{\}$ and $U_i \cap V = \{\}$.
c) In the sense of Gödel theorem it is not decidable that either one element of V exists or that $V = \{\}$.

Note:
(xiv) Criterion c) means that without V the equilibrium for W in b) only can be
$W \geq U_1 \cup U_2 \cup \ldots \cup U_n$

Conclusion:

The main idea for those virtual subset is to enlarge the set theory.
In general it is neither decidable if V has in minimum one element nor if V have finite numbers of elements or infinite many.
That leads to the possibility of a virtual set having finite elements but this elements are not countable. Those finite but not countable virtual subsets will be a new category of the set theory.

## Appendix

The virtual set **I** can be used to get another perspective of the unsolved conjecture about Wieferich prime numbers:

$2^{p-1} \equiv 1 \bmod p^2$

Only the both primes 1093 and 3511 are known.

$w_1 = 1093 \in$ **B** (brother is 1091).
$w_2 = 3511 \in$ **C** (cousin is $3511 + 2^{44}$).

With the virtual set **I** we can make the conjecture:

If there is a third Wieferich prime number $w_3$ then

$w_3 \in$ **I**

or

if **I** = { } then no more Wieferich prime numbers exist.

So if in future a third Wieferich number will be found, we can try to find out if this Wieferich number is a candidate of **I**.